\documentclass[12pt,a4paper]{amsart}
\usepackage{amsfonts,amssymb,amscd,amsmath,enumerate,verbatim,calc}
\usepackage{mathrsfs}
\usepackage[all]{xy}
\headsep=1cm \numberwithin{equation}{section}
\newtheorem{thm}{Theorem}[section]
\newtheorem{cor}[thm]{Corollary}
\newtheorem{lem}[thm]{Lemma}

\newtheorem{prop}[thm]{Proposition}
\theoremstyle{definition}

\theoremstyle{remark}

\numberwithin{equation}{section}

\newcommand\Supp{\operatorname{Supp}}
\newcommand\Ass{\operatorname{Ass}}
\newcommand\Assh{\operatorname{Assh}}
\newcommand\mAss{\operatorname{mAss}}

\newcommand\Ann{\operatorname{Ann}}

\newcommand\Spec{\operatorname{Spec}}
\newcommand\Rad{\operatorname{Rad}}
\newcommand\cd{\operatorname{cd}}
\newcommand\Hom{\operatorname{Hom}}
\newcommand\Ext{\operatorname{Ext}}
\newcommand\Tor{\operatorname{Tor}}

\newcommand\height{\operatorname{height}}
\newcommand\Att{\operatorname{Att}}
\newcommand\Max{\operatorname{Max}}

\newcommand\m{\operatorname{\frak m}}

\newcommand\p{\operatorname{\frak p}}
\newcommand\q{\operatorname{\frak q}}
\newcommand\Q{\operatorname{\frak Q}}

\begin{document}\title[A study of cofiniteness]{A study of cofiniteness through minimal associated primes}
\author[K. Bahmanpour]{  Kamal Bahmanpour }

\address{Department of Mathematics, Faculty of Sciences, University of Mohaghegh Ardabili,
56199-11367, Ardabil, Iran;
and School of Mathematics, Institute for Research in Fundamental Sciences (IPM), P.O. Box. 19395-5746, Tehran, Iran.} \email{\it bahmanpour.k@gmail.com}

\thanks{ This research of the author was supported by a grant from IPM (No. 96130018).}

\begin{abstract}
  In this paper we shall investigate the concepts of cofiniteness of
   local cohomology modules and Abelian categories of cofinite modules
   over arbitrary Noetherian rings. Then we shall improve some of the
   results given in \cite{B, BN, Ch, De, DM, DFT, Ha, HK, Ka0, Ka1, Ka2, Me2, PAB}.
  \end{abstract}
\maketitle
\section{Introduction}

Throughout this paper, let $R$ denote a commutative Noetherian ring
(with identity) and $I$ an ideal of $R$. For an $R$-module $M$, the
$i$-th local cohomology module of $M$ with support in $V(I)$
is defined as:
$$H^i_I(M) = \underset{n\geq1} {\varinjlim}\,\, \Ext^i_R(R/I^n,
M).$$   In fact, the local cohomology functors $H^i_I(-)$ arise as the
derived functors of the left exact
functor $\Gamma _I(-)$, where for an $R$-module $M$, $\Gamma _I(M)$
is the submodule of $M$ consisting of all elements annihilated by
some power of $I$, i.e., $\bigcup_{n=1}^{\infty}(0:_MI^n)$. Local
cohomology was first defined
and studied by Grothendieck. We refer the reader to \cite{BS} or \cite{Gr1}
 for more
details about local cohomology.\\

 For an $R$-module $M$, the notion $\cd (I, M)$, the
 cohomological dimension of $M$ with respect to $I$, is defined to be the
greatest integer $i$ such that $H^i_I(M)\neq0$ if there exist such $i'$s and $-\infty$
otherwise. Hartshorne \cite{Ha2} has defined the notion $q(I,R)$ as the
greatest integer $i$ such that $H^i_I(R)$
is not Artinian. Dibaei and Yassemi \cite{DY} extended this notion to
arbitrary $R$-modules, to
the effect that for any $R$-module $N$ they defined $q(I,N)$ as the greatest
integer $i$ such that
$H^i_I(N)$ is not Artinian if there exist such $i'$s and $-\infty$.
otherwise.\\

Hartshorne in \cite{Ha}  defined an $R$-module $L$ to be
$I$-{\it cofinite}, if $\Supp L\subseteq
V(I)$ and ${\rm Ext}^{i}_{R}(R/I, L)$ is a finitely generated module
for all $i$. Then he posed the following questions:

(i) {\it For which Noetherian rings $R$ and ideals $I$ are the modules
$H^{i}_{I}(M)$ $I$-cofinite for  all finitely generated
$R$-modules $M$ and all $i\geq0$}?

(ii) {\it Whether the category $\mathscr{C}(R, I)_{cof}$
of $I$-cofinite modules is an Abelian subcategory of the category of
 all $R$-modules?
That is, if $f: M\longrightarrow N$ is an $R$-homomorphism of
$I$-cofinite modules, are $\ker f$ and ${\rm coker} f$ $I$-cofinite?}\\

With respect to the question (i), Hartshorne in \cite{Ha} and later
Chiriacescu in \cite{Ch} showed that if $R$ is a complete regular
local ring and $I$ is a prime ideal such that $\dim R/I=1$, then
$H^{i}_{I}(M)$ is $I$-cofinite for any finitely generated $R$-module
$M$. This result was later extended to more general local rings and
one-dimensional ideals by Huneke and
Koh in \cite{HK} and by Delfino in \cite{De} until finally Delfino
and Marley in \cite{DM} and Yoshida
in \cite{Yo} proved that the local cohomology modules $H^i_I(M)$ are
 $I$-cofinite for all finitely generated $R$-modules $M$, where the
ideal $I$ of a local ring $R$, satisfies $\dim R/I = 1.$ Finally,
the local condition on the ring has been removed in \cite{BN}.
For some other related results, see also \cite{BNS0, MV}.
Furthermore, with respect to the question (i), Kawasaki in \cite{Ka0}
proved that  if an ideal $I$ of a Noetherian ring $R$ is principal,
up to radical, then the local cohomology modules $H^i_I(M)$ are
$I$-cofinite, for all finitely generated $R$-modules $M$ and
all integers $i\geq0$.  Also, Melkersson in \cite{Me} extended
this result for all ideals $I$ of $R$ with $\cd(I,R)\leq 1$.
Finally, the present author in \cite{B} generalized this
result for all ideals $I$ with $q(I,R)\leq 1$.\\

With respect to the question (ii), Hartshorne gave a counterexample
to show that this question has not an
affirmative answer in general, (see \cite[\S 3]{Ha}).
On the positive side, Hartshorne proved that if $I$ is a prime
ideal of dimension one in a complete regular local ring $R$,
then the answer to his question
is yes. On the other hand, in \cite{DM}, Delfino and Marley
extended this result to arbitrary
Noetherian complete local rings. Kawasaki in \cite{Ka2} generalized
the Delfino and Marley result
for an arbitrary ideal $I$ of dimension one in a local ring $R$. Finally, Melkersson
in \cite{Me2} generalized the Kawasaki's result for all ideals of dimension one of any
arbitrary Noetherian ring $R$. Furthermore, in \cite{BNS} as a generalization of
Melkersson's result it is shown that for any ideal $I$ in a Noetherian ring $R$,
the category of all $I$-cofinite
$R$-modules $M$ with $\dim M\leq 1$ is an Abelian subcategory of the category of all
$R$-modules. For some other similar results,  see also \cite{B1}. Moreover, with
respect to the question (ii),
Kawasaki in \cite{Ka1}
proved that if an ideal $I$ of a Noetherian ring $R$ is principal, up to radical,
then the category $\mathscr{C}(R, I)_{cof}$ is Abelian. Pirmohammadi et al in \cite{PAB}
as a generalization of Kawasaki's result proved that if $I$ is an ideal of a
Noetherian local ring $R$ with $\cd(I,R)\leq1$, then $\mathscr{C}(R, I)_{cof}$
is Abelian. Also, more recently, Divaani-Aazar et al. in \cite{DFT} have removed
the local condition on the ring. Also, the present author in \cite{B} proved that if $I$ is
an ideal of a  Noetherian complete local ring $R$ with $q(I,R)\leq1$, then $\mathscr{C}(R, I)_{cof}$ is Abelian.\\

Now, for any ideal $I$ of $R$ and any finitely generated $R$-module $M$ we define $$\mathfrak{A}(I,M):=\{\p \in \mAss_R M\,\,:\,\,I+\p=R\,\,\,\,{\rm or}\,\,\,\,\p\supseteq I\},$$$$\mathfrak{B}(I,M):=\{\p \in \mAss_R M\,\,:\,\,\cd(I,R/\p)=1\},$$$$\mathfrak{C}(I,M):=\{\p \in \mAss_R M\,\,:\,\,q(I,R/\p)=1\}\,\,{\rm and}$$$$\mathfrak{D}(I,M):=\{\p \in \mAss_R M\,\,:\,\,0\leq \dim R/(I+\p)\leq 1\}.$$

Also, we denote by $\mathscr{I}(R)$ the class of all ideals $I$ of $R$ with the property that, for every finitely generated $R$-module $M$, the local cohomology modules $H^i_I(M)$ are $I$-cofinite for all $i\geq 0$. \\

In Section 2 of this paper as a generalization of the some results given in \cite{DM, DFT, Ha, Ka0, Ka1, Ka2, Me2, PAB} we shall prove the following result:\\

{\bf Theorem 1}. {\it  Let $I$ be an ideal of $R$ such that $$\mAss_R R=\mathfrak{A}(I,R)\cup\mathfrak{B}(I,R)\cup\mathfrak{D}(I,R).$$  Then the category
$\mathscr{C}(R, I)_{cof}$ is Abelian.}\\

In Section 3 we present a generalization of the some results given in \cite{B, BN, Ch, De, DM, Ha, HK, Yo} as follows:\\

 {\bf Theorem 2}. {\it  Let $I$ be an ideal of $R$ such that $$\mAss_R R=\mathfrak{A}(I,R)\cup\mathfrak{B}(I,R)\cup \mathfrak{C}(I,R)\cup\mathfrak{D}(I,R).$$ Then $I\in \mathscr{I}(R)$.}\\

In Section 4, we shall present a formula for the cohomological dimension of finitely generated modules with respect to ideals of a Noetherian complete local ring $R$ belong to $\mathscr{I}(R)$ in terms of the height of ideals. More precisely, we shall prove the following result:\\

{\bf Theorem 3}. {\it  Let $(R,\m)$ be a Noetherian complete local ring, $I\in \mathscr{I}(R)$ and $M$ be a non-zero finitely generated $R$-module.  Then,
  $$\cd(I,M)={\rm max}\{\height (I+\p)/\p\,\,:\,\,\p\in \mAss_R M\}.$$}\\

Throughout this paper, for any ideal $I$ of $R$, we denote by $\mathscr{C}(R, I)_{cof}$ the category of all $I$-cofinite
$R$-modules. Also, for each $R$-module $L$, we denote by
 $\Assh_R L$ (respectively by $\mAss_R L$), the set $\{\p\in \Ass_R L\,\,:\,\, \dim R/\p= \dim L\}$  (respectively the set of
 minimal elements of $\Ass_R L$ with respect to inclusion). For any ideal $\frak{a}$ of $R$, we denote
$\{\frak p \in {\rm Spec}\,R:\, \frak p\supseteq \frak{a}\}$ by
$V(\frak{a})$. For any ideal $\frak{b}$ of $R$, {\it the
radical of} $\frak{b}$, denoted by $\Rad(\frak{b})$, is defined to
be the set $\{x\in R \,: \, x^n \in \frak{b}$ for some $n \in
\mathbb{N}\}$.  Finally, we denote by $\Max(R)$ the set of all maximal ideals of $R$. Furthermore, in this paper we interpret the Krull dimension of the zero module as $-\infty$.  For any unexplained notation and terminology we refer the reader to \cite{BS} and \cite{Mat}.

\section{Abelianness of the category of cofinite modules}

In this section we extend some results given in  \cite{DM, DFT,  Ha, Ka0, Ka1, Ka2, Me2, PAB}. The main purpose of this section is to prove Theorem 2.5. The following lemma and proposition are needed for the proof of Corollary 2.3. We recall that for any proper ideal $I$ of $R$, the {\it arithmetic rank } of $I$,
denoted by ${\rm ara}(I)$, is the least number of elements of $I$
required to generate an ideal which has the same radical as $I$.

 \begin{lem}
 \label{2.1}
 Let $I$ be an ideal of $R$ and $M$ be an $I$-cofinite
$R$-module. Then for each finitely generated $R$-module $N$ with $\Supp N/IN \subseteq \Max(R)$
the $R$-modules $\Ext^i_R(N,M)$ and $\Tor^R_i(N,M)$ are Artinian and $I$-cofinite, for all
$i\geq 0$.
  \end{lem}
\proof By the similarity of the proof we prove the assertion just for the $R$-modules $\Ext^i_R(N,M)$, $i\geq 0$.  If $IN=N$ then for each $i\geq 0$ $$\emptyset=\Supp N/IN=\Supp N\cap V(I)\supseteq \Supp N\cap \Supp M\supseteq \Supp \Ext^i_R(N,M).$$ Hence, in this case we have $\Ext^i_R(N,M)=0$ for all integers $i\geq0$ and so the assertion is clear. So, without loss of generality we may assume that $N/IN\neq0$. Then we have $I+\Ann_R N\neq R$. Now in order to prove the assertion we use induction on $$t={\rm ara}(I+{\rm Ann}_R N/{\rm
Ann}_R N).$$ If $t=0$ then it follows from the definition that
${\rm Supp}(N)\subseteq V(I)$ and so it follows from \cite[Corollary 1]{DM} or \cite[Corollary 2.5]{Me} that, for each integer $i\geq0$, the $R$-module $\Ext^i_R(N,M)$ is finitely generated with support in the set  $\Supp N\cap \Supp M$. But it is clear that $$ \Supp N\cap \Supp M\subseteq \Supp N\cap V(I)=\Supp N/IN\subseteq \Max(R)$$ which means that for each integer $i\geq0$, the $R$-module $\Ext^i_R(N,M)$ is of finite length. So assume that $t>0$ and the result has been proved for
$0,1,...,t-1$. Since ${\rm Ann}_R N\subseteq {\rm
Ann}_R N/\Gamma_{I}(N)$, it follows that $${\rm ara}(I+{\rm
Ann}_R N/\Gamma_{I}(N)/{\rm Ann}_R N/\Gamma_{I}(N))\leq {\rm
ara}(I+{\rm Ann}_R N/{\rm Ann}_R N).$$On the other hand, the exact
sequence
$$0\longrightarrow \Gamma_{I}(N) \longrightarrow N \longrightarrow N/\Gamma_{I}(N)
\longrightarrow 0,$$ induces the following exact
sequence$$0\longrightarrow {\rm
Hom}_R(N/\Gamma_{I}(N),M)\longrightarrow {\rm
Hom}_R(N,M)\longrightarrow {\rm
Hom}_R(\Gamma_{I}(N),M)$$$$\longrightarrow{\rm
Ext}^1_R(N/\Gamma_{I}(N),M)\longrightarrow{\rm
Ext}^1_R(N,M)\longrightarrow {\rm
Ext}^1_R(\Gamma_{I}(N),M)\longrightarrow \cdots.\,\,\,\,(2.1.1)$$But, using the facts that $$\Supp \Gamma_{I}(N)/I \Gamma_{I}(N)=\Supp \Gamma_{I}(N)\cap V(I)\subseteq \Supp N\cap V(I)=\Supp N/IN\subseteq \Max(R)$$ and
${\rm ara}(I+{\rm Ann}_R \Gamma_{I}(N)/{\rm
Ann}_R \Gamma_{I}(N))= 0$, the inductive hypothesis  yields that the $R$-modules ${\rm Ext}^i_R(\Gamma_{I}(N),M)$ are of finite length and $I$-cofinite, for all integers $i\geq 0$. So, using the exact sequence $(2.1.1)$, \cite{BNS} and by replacing $N$ by $N/\Gamma_{I}(N)$, without loss of generality, we may
assume that $N$ is a finitely generated $I$-torsion-free $R$-module with $$\emptyset\neq\Supp N/IN\subseteq \Max(R)\,\,{\rm and}\,\,{\rm ara}(I+{\rm
Ann}_R N/{\rm Ann}_R N)=t.$$
Then, by the definition there exist elements
$y_1,...,y_t\in I$, such that $${\rm Rad}(I+{\rm Ann}_R N/{\rm
Ann}_R N)={\rm Rad}((y_1,...,y_t)+{\rm Ann}_R N/{\rm Ann}_R N).$$
Furthermore, by \cite [Lemma 2.1.1]{BS}, $I
\nsubseteq \bigcup_{\frak p \in \Ass_{R} N} \frak p$.

Therefore,

\begin{center}
$(y_1,...,y_t)+{\rm Ann}_R N \not \subseteq \bigcup_{\frak p \in
\Ass_{R}N} \frak p.$
\end{center}

But, as

\begin{center}
${\rm Ann}_R N \subseteq\bigcap_{\frak p \in \Ass_{R}N} \frak p,$
\end{center}

it follows that

\begin{center}
$(y_1,...,y_t)\not \subseteq \bigcup_{\frak p \in \Ass_{R}N} \frak p.$
\end{center}

Therefore, by \cite [Exercise 16.8]{Mat} there is $a\in (y_2,\dots,y_t)$
such that
\begin{center}
$y_1+a\not\in\bigcup_{\frak p \in \Ass_{R}N} \frak p.$
\end{center}
Let $x:=y_1+a$. Then $x\in I$ and $$\Rad(I+{\rm Ann}_R N/{\rm
Ann}_R N)={\rm Rad}((x,y_2,...,y_t)+{\rm Ann}_R N/{\rm
Ann}_RN).$$ Now it is easy to see that $$\Rad(I+{\rm
Ann}_R N/xN/{\rm Ann}_R N/xN)={\rm Rad}((y_2,...,y_t)+{\rm
Ann}_R N/xN/{\rm Ann}_R N/xN).$$ and hence ${\rm ara}(I+{\rm
Ann}_R N/xN/{\rm Ann}_R N/xN)\leq t-1$.
Also, it is clear that $$\Supp (N/xN)/I (N/xN)=\Supp N/xN\cap V(I)\subseteq \Supp N\cap V(I)=\Supp N/IN\subseteq \Max(R).$$ Therefore, by the inductive hypothesis the $R$-module $\Ext^i_R(N/xN,M)$ is Artinian and $I$-cofinite for each integer $i\geq 0$.
 Now, the exact sequence
$$0\longrightarrow N  \stackrel{x} \longrightarrow N \longrightarrow
N/xN \longrightarrow 0$$ induces an  exact sequence $${\rm
Ext}^i_R(N,M) \stackrel{x}\longrightarrow {\rm Ext}^i_R(N,M)
\longrightarrow {\rm Ext}^{i+1}_R(N/xN,M)$$$$\longrightarrow  {\rm
Ext}^{i+1}_R(N,M) \stackrel{x}\longrightarrow {\rm
Ext}^{i+1}_R(N,M),$$ for all integers $i\geq 0$. Consequently, for all integers $i\geq 0$, we have the following short exact sequence,$$0
\longrightarrow {\rm Ext}^i_R(N,M)/x{\rm Ext}^i_R(N,M)
\longrightarrow {\rm Ext}^{i+1}_R(N/xN,M) \longrightarrow (0: _{{\rm
Ext}^{i+1}_R(N,M)} x) \longrightarrow 0.$$

But, the $R$-modules ${\rm
Ext}^{i+1}_R(N/xN,M),$ are $I$-cofinite and Artinian, for all $i\geq 0$ and
hence it follows from \cite[Corollary 4.4]{Me} that for all integers $i\geq0$, the $R$-modules ${\rm Ext}^i_R(N,M)/x{\rm Ext}^i_R(N,M)$ and $(0 :_{{\rm
Ext}^{i+1}_R(N,M)}x)$ are $I$-cofinite. Furthermore, it
follows from the exact sequence $$0\longrightarrow {\rm
Hom}_R(N/xN,M) \longrightarrow {\rm Hom}_R(N,M)\stackrel{x}
\longrightarrow {\rm Hom}_R(N,M),$$ and inductive hypothesis that
the $R$-module $(0 :_{{\rm Hom}_R(N,M)}x)$ is also $I$-cofinite.
 Now,  since the
$R$-modules $(0:_ {{\rm Ext}^i_R(N,M)} x)$ and ${\rm
Ext}^i_R(N,M)/x{\rm Ext}^i_R(N,M)$ are $I$-cofinite for all $i\geq 0$, it follows from \cite[Corollary 3.4]{Me} that the $R$-modules ${\rm
Ext}^i_R(N,M)$ are $I$-cofinite for all integers $i\geq0$. Moreover, since for each integer $i\geq0$, $$ \Supp \Ext^i_R(N,M)\subseteq \Supp N\cap \Supp M\subseteq \Supp N\cap V(I)=\Supp N/IN\subseteq \Max(R),$$ it follows that the $R$-module $$(0:_{\Ext^i_R(N,M)}I)\simeq \Hom_R(R/I,\Ext^i_R(N,M))$$ is of finite length and so the result \cite[Proposition 4.1]{Me} yields that the $R$-module $\Ext^i_R(N,M)$ is Artinian and $I$-cofinite, for
each integer $i\geq0$. This completes the inductive step and the proof of the lemma.\qed\\

\begin{prop}
 \label{2.2}
 Let $I$ be an ideal of $R$ and $M$ be an $I$-cofinite
$R$-module. Then for each finitely generated $R$-module $N$ with $\dim N/IN \leq 1$,
the $R$-modules $\Ext^i_R(N,M)$ and $\Tor^R_i(N,M)$ are $I$-cofinite, for all
$i\geq 0$.
\end{prop}

\proof By the similarity of the proof we prove the assertion just for the $R$-modules $\Ext^i_R(N,M)$, $i\geq 0$.  We use induction on $t={\rm ara}(I+{\rm Ann}_R N/{\rm
Ann}_RN)$. If $t=0$, then, it follows from the definition that
${\rm Supp}(N)\subseteq V(I)$ and so the assertion holds by
 \cite[Corollary 2.5]{Me}. So assume that $t>0$ and the result has been proved for
$0,1,...,t-1$. Since ${\rm Ann}_RN\subseteq {\rm
Ann}_R N/\Gamma_{I}(N)$, it follows that $${\rm ara}(I+{\rm
Ann}_R N/\Gamma_{I}(N)/{\rm Ann}_R N/\Gamma_{I}(N))\leq {\rm
ara}(I+{\rm Ann}_R N/{\rm Ann}_RN).$$On the other hand, the exact
sequence
$$0\longrightarrow \Gamma_{I}(N) \longrightarrow N \longrightarrow N/\Gamma_{I}(N)
\longrightarrow 0,$$ induces the following exact
sequence$$0\longrightarrow {\rm
Hom}_R(N/\Gamma_{I}(N),M)\longrightarrow {\rm
Hom}_R(N,M)\longrightarrow {\rm
Hom}_R(\Gamma_{I}(N),M)$$$$\longrightarrow{\rm
Ext}^1_R(N/\Gamma_{I}(N),M)\longrightarrow{\rm
Ext}^1_R(N,M)\longrightarrow {\rm
Ext}^1_R(\Gamma_{I}(N),M)\longrightarrow \cdots.$$ So, using Lemma 2.1, \cite{BNS} and \cite[Corollary 2.5]{Me}, by replacing $N$ by $N/\Gamma_{I}(N)$,  without loss generality, we may
assume that $N$ is a finitely
generated $I$-torsion-free $R$-module, such that $\dim N/IN=1$ and ${\rm ara}(I+{\rm
Ann}_R N/{\rm Ann}_R N)=t$. Then, by \cite [Lemma 2.1.1]{BS}, $I
\nsubseteq \bigcup_{\frak p \in \Ass_{R} N} \frak p$. Next, let
$k\geq 0$ and
\begin{center}
$S_k:=\bigcup_{i=0}^{k}\Supp {\rm Ext}^i_R(N,M),$
\end{center}

and $$T:=\{\frak p\in S_k \mid \dim R/\frak p= 1\}.$$  Now, it is
easy to see that $T\subseteq {\rm Assh}_R N/IN$. Therefore $T$ is a finite set. Moreover,
for each ${\frak p}\in T$, using \cite[Exerxise 7.7]{Mat} it follows that $M_{\frak p}$ is a $IR_{\frak p}$-cofinite
module and $N_{\p}$ is a finitely generated $R_{\p}$-module with $$\Supp N_{\p}/IN_{\p}=\{\p R_{\p}\}=\Max(R_{\p}).$$ Therefore, using \cite[Exerxise 7.7]{Mat} and Lemma 2.1 it follows that the $R_{\frak p}$-module $({\rm Ext}^i_R(N,M))_{\frak p}$ is Artinian and $IR_{\frak
p}$-cofinite, for each $i\geq0$.  Now, applying the method used in the proof of \cite[Theorem 2.7]{AB} with the same notation it follows that the $R$-modules ${\rm
Ext}^i_R(N,M)$ are $I$-cofinite for all $0\leq i \leq k$. Therefore,
as $k$ is arbitrary, it follows that, the $R$-modules ${\rm
Ext}^i_R(N,M)$ are of $I$-cofinite, for all $i\geq 0$.
This completes the inductive step and the proof of theorem. \qed\\

We need the following consequence of Proposition 2.2 in the proof of Theorem 2.5.

\begin{cor}
 \label{2.3}
 Let $I$ be an ideal of $R$ and $M$ be an $I$-cofinite
$R$-module. Set $J:= \bigcap _{\p\in \mathfrak{D}(I,R)} \p$ and $K:=\Gamma_J(R)$. Then, the $R$-modules $\Ext^i_R(K,M)$ are $I$-cofinite for all $i\geq 0$.
\end{cor}
\proof Since, $\Ass_R K=\Ass_R R \cap V(J)$ it follows that $\mAss_R K = \mathfrak{D}(I,R)=\mathfrak{D}(I,K)$. Therefore, for each $\p\in \mAss_R K$ by the definition we have $0\leq\dim R/(I+\p)\leq 1$. Thus, we have $0\leq\dim K/IK\leq 1$. So, the assertion follows from Proposition 2.2.\qed\\

The following lemma is needed in the proof of Theorem 2.5.

\begin{lem}
 \label{2.4}
 Let $I$ and $J$ be two proper ideals of $R$ and $M$ be an $R$-module with $JM=0$ and $\Supp M \subseteq V(I)$. Then $M$ as an $R$-module is $I$-cofinite if and only if $M$ as an $R/J$-module is $(I+J)/J$-cofinite.
 \end{lem}
 \proof See \cite[Proposition 2]{DM}.\qed\\

Now we are ready to state and prove the main result of this section.

\begin{thm}
 \label{2.5}
Let $I$ be an ideal of $R$ such that $$\mAss_R R=\mathfrak{A}(I,R)\cup\mathfrak{B}(I,R)\cup\mathfrak{D}(I,R).$$  Then
$\mathscr{C}(R, I)_{cof}$ is Abelian.
\end{thm}
\proof Let $M,N\in \mathscr{C}(R, I)_{cof}$ and let $f:M\longrightarrow N$ be
an $R$-homomorphism. It is enough to prove that the $R$-modules ${\rm ker}\,f$ and ${\rm coker}\,f$
are $I$-cofinite.\\

We consider the following three cases:\\

Case 1. Assume that $\mathfrak{D}(I,R)=\emptyset$. Then $T:=\bigoplus_{\p\in \mAss_R R}R/\p$ is a finitely generated
 $R$-module with $\Supp T = \Spec R = \Supp R$. So, it follows from \cite[Theorem 2.2]{DNT}
  that $$\cd(I,R)=\cd(I,T)={\rm sup}\{\cd(I,R/\p)\,\,:\,\,\p\in \mathfrak{A}(I,R)\cup\mathfrak{B}(I,R)\}\leq 1.$$ Thus, the assertion follows from \cite[Theorem 2.2]{DFT}.\\

Case 2. Assume that $\mathfrak{A}(I,R)\cup\mathfrak{B}(I,R)=\emptyset$. Then since by the hypothesis for each $\p\in \mathfrak{D}(I,R)$ we have
 $0\leq\dim R/(I+\p)\leq 1$ it follows that for each
$\p\in \mAss_R R$ we have $0\leq\dim R/(I+\p)\leq 1$. Thus $\dim(R/I)\leq 1$ and hence the assertion follows from \cite[Theorem 2.7]{BNS}. \\

Case 3. Assume that $\mathfrak{A}(I,R)\cup\mathfrak{B}(I,R)\neq\emptyset$ and $\mathfrak{D}(I,R)\neq\emptyset$. Set $$J:=\bigcap_{\p\in \mathfrak{D}(I,R)}\p\,\,{\rm and}\,\,K:=\Gamma_{J}(R).$$ Then by Corollary 2.3 the $R$-modules $\Hom_R(K,M)$ and $\Hom_R(K,N)$ are $I$-cofinite. The exact sequence $$0 \longrightarrow K \stackrel{\iota}\longrightarrow R \stackrel{\pi}\longrightarrow R/K \longrightarrow0,$$ induces the following commutative diagrams with exact rows
\begin{displaymath}
\xymatrix{0\ar[r]&\Hom_R(R/K,M)\ar[r]\ar[d]& \Hom_R(R,M)\ar[d]\ar[r]&\Hom_R(K,M)\ar[d]\\
0\ar[r]&\Hom_R(R/K,N)\ar[r]& \Hom_R(R,N) \ar[r]&\Hom_R(K,N).&\,\,\,\,\,\,\,\,(2.5.1) }
\end{displaymath}
and
\begin{displaymath}
\xymatrix{\Hom_R(R,M)\ar[r]\ar[d]& \Hom_R(K,M)\ar[d]\ar[r]&\Ext^1_R(R/K,M)\ar[r]\ar[d]&0\\
\Hom_R(R,N)\ar[r]& \Hom_R(K,N) \ar[r]&\Ext^1_R(R/K,N)\ar[r]&0&\,\,\,\,\,\,\,\,(2.5.2) }
\end{displaymath}

Two diagrams $(2.5.1)$ and $(2.5.2)$ induce the exact sequences
$$0\longrightarrow {\rm im }\,\Hom_R (\iota,M)\longrightarrow \Hom_R(K,M) \longrightarrow  \Ext^1_R(R/K,M)\longrightarrow 0,\,\,\,\,\,(2.5.3)$$and
$$0\longrightarrow \Hom_R(R/K,M)\longrightarrow \Hom_R(R,M) \longrightarrow  {\rm im }\,\Hom_R (\iota,M)\longrightarrow 0.\,\,\,\,\,(2.5.4)$$Also, the diagram $(2.5.1)$ induces a commutative diagram with exact rows
\begin{displaymath}
\xymatrix{0\ar[r]&\Hom_R(R/K,M)\ar[r]\ar[d]& \Hom_R(R,M)\ar[d]\ar[r]&{\rm im }\,\Hom_R (\iota,M)\ar[d]^{\alpha}\ar[r]&0\\
0\ar[r]&\Hom_R(R/K,N)\ar[r]& \Hom_R(R,N) \ar[r]&{\rm im }\,\Hom_R (\iota,N)\ar[r]&0.&\,\,\,\,\,\,\,\,(2.5.5) }
\end{displaymath}
Set $T:=\Ext^1_R(R/K,M)$. The exact sequence
$(2.5.3)$
induces an  exact sequence

$$\Tor^R_1(R/(I+K),\Hom_R(K,M))\longrightarrow \Tor^R_1(R/(I+K),T)$$$$\longrightarrow \Tor^R_0(R/(I+K),{\rm im }\,\Hom_R (\iota,M))
\rightarrow  \Tor^R_0(R/(I+K),\Hom_R(K,M))$$$$\longrightarrow  \Tor^R_0(R/(I+K),T)\longrightarrow 0.\,\,\,\,\,(2.5.6)$$ By \cite[Theorem 2.1 and Corollary 2.5]{Me}
 the modules $$\Tor^R_0(R/(I+K),M),\,\,\Tor^R_1(R/(I+K),\Hom_R(K,M))\,\,{\rm and}\,\,\Tor^R_0(R/(I+K),\Hom_R(K,M))$$ are finitely generated $R$-modules. Furthermore, the exact sequence $$M  \longrightarrow {\rm im }\,\Hom_R (\iota,M)\longrightarrow 0 $$induces an exact sequence $$\Tor^R_0(R/(I+K),M)\longrightarrow \Tor^R_0(R/(I+K),{\rm im }\,\Hom_R (\iota,M))\longrightarrow 0,$$ which implies the $R$-module $\Tor^R_0(R/(I+K),{\rm im }\,\Hom_R (\iota,M))$ is finitely generated. So, from the exact sequence $(2.5.6)$ we deduce that the $R$-modules $$\Tor^R_0(R/(I+K),T)\,\,{\rm and}\,\,\Tor^R_1(R/(I+K),T)$$ are finitely generated. By the same method it follows that the $R$-modules $$\Tor^R_0(R/I,T)\,\,{\rm and}\,\,\Tor^R_1(R/I,T)$$ are finitely generated.\\

   Since the $R$-module $T/IT\simeq \Tor^R_0(R/I,T)$ is finitely generated it follows that there is a finitely generated submodule $T_1$ of $T$ such that $T/IT=(T_1+IT)/IT$ and hence $T=T_1+IT$. Since, $T_1$ is an $I$-torsion finitely generated $R$-module it follows that $I^kT_1=0$ for some positive integer $k$. Therefore, $I^kT=I^kT_1+I^{k+1}T=I^{k+1}T$. Set $L:=I^kT$. Then $IL=L$ and $T/L$ is finitely generated; because the $R$-module $T/IT$ is finitely generated.\\

     The exact sequence $$0 \longrightarrow L \longrightarrow T \longrightarrow T/L \longrightarrow0\,\,\,\,(2.5.7)$$ induces an exact sequence $$\Tor^R_2(R/(I+K),T/L)\longrightarrow\Tor^R_1(R/(I+K),L)\longrightarrow\Tor^R_1(R/(I+K),T)$$ and hence the $R$-module $\Tor^R_1(R/(I+K),L)$ is finitely generated. Also, applying the same method it follows that the $R$-module $\Tor^R_1(R/I,L)$ is finitely generated. Since, the $R$-module $(I\cap K)/IK$ has an $R/I$-module structure it follows that for some positive integer $n$ there is an exact sequence $$\bigoplus^n_{i=1}R/I\longrightarrow (I\cap K)/IK \longrightarrow0,$$which yields an exact sequence $$\bigoplus^n_{i=1}L/IL\longrightarrow (I\cap K)/IK\otimes_R L\longrightarrow0.$$ But, we have $L/IL=0$ and so $(I\cap K)/IK\otimes_R L=0$. Furthermore, the exact sequence $$0\longrightarrow (I\cap K)/IK\longrightarrow I/IK \longrightarrow I/(I\cap K)\longrightarrow 0$$ yields the exact sequence $$(I\cap K)/IK\otimes_R L\longrightarrow I/IK \otimes_R L\longrightarrow I/(I\cap K)\otimes_R L\longrightarrow 0,$$whence we get the isomorphisms $$ (I+K)/K \otimes_R L \simeq I/(I\cap K)\otimes_R L \simeq I/IK \otimes_R L.$$
      Using the fact that $KL=0$ we have $L\simeq L/KL\simeq R/K\otimes_R L$. Therefore, $$I\otimes_R L\simeq I\otimes_R(R/K\otimes_RL )\simeq (I\otimes_RR/K)\otimes_RL\simeq I/KI \otimes_RL.$$
 Therefore, $$I\otimes_R L\simeq (I+K)/K \otimes_R L.$$
The exact sequence $$0 \longrightarrow I\longrightarrow R \longrightarrow R/I \longrightarrow0$$induces the exact sequence $$0 \longrightarrow \Tor^R_1(R/I,L)\longrightarrow I\otimes_R L \longrightarrow L\longrightarrow0,$$
which using the fact that $L=IL=(I+K)L$ induces an exact sequence $$\Tor^R_1(R/I,L)\otimes_R R/(I+K)\longrightarrow (I\otimes_R L)\otimes_R R/(I+K)\longrightarrow0, $$ whence we conclude that the $R$-module $$((I+K)/K \otimes_R L)\otimes_R R/(I+K) \simeq (I\otimes_R L) \otimes_R R/(I+K) $$is finitely generated. Moreover, the exact sequence $$0 \longrightarrow (I+K)/K\longrightarrow R/K \longrightarrow R/(I+K) \longrightarrow0$$induces the exact sequence $$0 \longrightarrow \Tor^{R/K}_1(R/(I+K),L)\longrightarrow (I+K)/K\otimes_{R/K} L \longrightarrow L\longrightarrow0.$$The last exact sequence using the facts that $L=IL=(I+K)L$ and $(I+K)/K\otimes_{R/K} L\simeq (I+K)/K\otimes_{R} L$
 induces the exact sequence $$\Tor^{R}_1(R/(I+K),L)\longrightarrow \Tor^{R/K}_1(R/(I+K),L)\longrightarrow ((I+K)/K\otimes_R L) \otimes_R R/(I+K)\longrightarrow0, $$
 which implies that the $R$-module $$\Tor^{R/K}_1(R/(I+K),L) $$is finitely generated. Furthermore, the exact sequence $(2.5.7)$ induces an exact sequence $$\Tor^{R/K}_2(R/(I+K),T/L)\longrightarrow\Tor^{R/K}_1(R/(I+K),L)\longrightarrow\Tor^{R/K}_1(R/(I+K),T)$$$$
 \longrightarrow\Tor^{R/K}_1(R/(I+K),T/L)\longrightarrow\Tor^{R/K}_0(R/(I+K),L)\longrightarrow$$
 $$\Tor^{R/K}_0(R/(I+K),T)\longrightarrow\Tor^{R/K}_0(R/(I+K),T/L)\longrightarrow0,$$
which implies that the $S$-modules $\Tor^{S}_0(S/IS,T)$ and $\Tor^{S}_1(S/IS,T)$ are finitely generated, where $S:=R/K$. But, since $\Ass_R\, R/K=\Ass_R\, R\backslash V(J)$ it follows that
 $\mAss_S\, S\subseteq\{\p S\,\,:\,\,\p\in \mathfrak{A}(I,R)\cup\mathfrak{B}(I,R)\}$. Whence, we can deduce that $\cd(IS,S)\leq 1$. Now, it follows from \cite[Theorem 2.9]{DFT2} that the $S$-module $T= \Ext^1_R(R/K,M)$ is $IS$-cofinite. Therefore, in view of Lemma 2.4, the $R$-module $\Ext^1_R(R/K,M)$ is $I$-cofinite. Now, it follows from the exact sequence $(2.5.3)$ that the $R$-module ${\rm im }\,\Hom_R (\iota,M)$ is $I$-cofinite. Also, it follows from the exact sequence $(2.5.4)$ that the $R$-module $\Hom_R(R/K,M)$ is $I$-cofinite. In particular, the $S$-module $\Hom_R(R/K,M)$ is $IS$-cofinite, by Lemma 2.4. By the same method we can prove that the $R$-modules $\Hom_R(R/K,N)$, $\Ext^1_R(R/K,N)$ and  ${\rm im }\,\Hom_R (\iota,N)$ are $I$-cofinite. In particular, the $S$-modules $\Hom_R(R/K,N)$ and $\Ext^1_R(R/K,N)$ are $IS$-cofinite, by Lemma 2.4. Since, the $S$-modules $\Hom_R(R/K,M)$ and $\Hom_R(R/K,N)$ are $IS$-cofinite and $\cd(IS,S)\leq 1$, it follows from   \cite[Theorem 2.2]{DFT} that the $S$-module ${\rm ker}\,\Hom_R(R/K,f)$ and ${\rm coker}\,\Hom_R(R/K,f)$ are $IS$-cofinite. In partiular, the $R$-modules ${\rm ker}\,\Hom_R(R/K,f)$ and ${\rm coker}\,\Hom_R(R/K,f)$ are $I$-cofinite, by Lemma 2.4. Also, we have $\dim {\rm im }\,\Hom_R (\iota,M)\leq 1 $, $\dim {\rm im }\,\Hom_R (\iota,N)\leq 1 $ and these $R$-modules are $I$-cofinite; so the result \cite[Theorem 2.7]{BNS} implies that the $R$-modules ${\rm ker}\, \alpha$ and ${\rm coker}\, \alpha$ are $I$-cofinite.\\

    Applying the Snake Lemma to the diagram $(2.5.5)$ we get an exact sequence $$0\longrightarrow{\rm ker}\,\Hom_R(R/K,f) \longrightarrow {\rm ker}\,f\stackrel{\mu}\longrightarrow {\rm ker}\,\alpha\stackrel{\lambda}\longrightarrow
{\rm coker}\,\Hom_R(R/K,f).\,\,\,\,(2.5.8)$$
The exact sequence $$0\longrightarrow {\rm im}\,\mu \longrightarrow {\rm ker}\,\alpha \longrightarrow {\rm im}\,\lambda \longrightarrow 0, $$
induces the exact sequence $$0\longrightarrow \Hom_R(R/I,{\rm im}\,\mu) \longrightarrow \Hom_R(R/I,{\rm ker}\,\alpha) \longrightarrow$$$$ \Hom_R(R/I,{\rm im}\,\lambda ) \longrightarrow \Ext^1_R(R/I,{\rm im}\,\mu) \longrightarrow  \Ext^1_R(R/I,{\rm ker}\,\alpha),$$which implies that the $R$-modules $\Hom_R(R/I,{\rm im}\,\mu)$ and $\Ext^1_R(R/I,{\rm im}\,\mu)$ are finitely generated. Since, $$\dim {\rm im}\,\mu \leq \dim {\rm ker}\,\alpha\leq \dim \Hom_R(K,M)\leq 1$$and the $R$-module ${\rm im}\,\mu$ is $I$-torsion it follows from \cite[Proposition 2.6]{BNS} that the $R$-module ${\rm im}\,\mu$ is $I$-cofinite. Moreover, the exact sequence $(2.5.8)$ yields an exact sequence $$0\longrightarrow{\rm ker}\,\Hom_R(R/K,f) \longrightarrow {\rm ker}\,f\longrightarrow {\rm im}\,\mu \longrightarrow0$$which implies that the $R$-module ${\rm ker}\,f$ is $I$-cofinite. Now, the exact sequences $$0\longrightarrow {\rm ker} \,f \longrightarrow M \longrightarrow
{\rm im}\, f \longrightarrow 0,$$ and
$$0\longrightarrow {\rm im}\,f \longrightarrow N \longrightarrow
{\rm coker}\,f \longrightarrow 0,$$ imply that the $R$-module ${\rm coker}\,f$ is $I$-cofinite too.\qed\\

\begin{cor}
 \label{2.6}
Let $I$ be an ideal of $R$ such that
$$\mAss_R R=\mathfrak{A}(I,R)\cup\mathfrak{B}(I,R)\cup\mathfrak{D}(I,R).$$
Let $$X^\bullet:
\cdots\longrightarrow X^i \stackrel{f^i} \longrightarrow X^{i+1}
\stackrel{f^{i+1}} \longrightarrow X^{i+2}\longrightarrow \cdots,$$ be a
complex such that  $X^i\in\mathscr{C}(R, I)_{cof}$ for all $i\in\Bbb{Z}$.
Then for each $i\in\Bbb{Z}$ the $i^{{\rm th}}$ cohomology module $H^i(X^\bullet)$ is in
$\mathscr{C}(R, I)_{cof}$.
\end{cor}
\proof The assertion follows from  Theorem 2.5.\qed\\

\begin{cor}
\label{2.7}
Let $I$ be an ideal of $R$ such that
$$\mAss_R R=\mathfrak{A}(I,R)\cup\mathfrak{B}(I,R)\cup\mathfrak{D}(I,R)$$ and
 $M$ be an $I$-cofinite $R$-module. Then,  the $R$-modules
$\Tor_i^R(N,M)$ and $\Ext^i_R(N,M)$ are $I$-cofinite, for all finitely generated $R$-modules $N$ and all
integers $i\geq0$.
\end{cor}
\proof Since $N$ is finitely generated it follows that, $N$ has a
free resolution with finitely generated free $R$-modules. Now the
assertion follows using Corollary 2.6 and computing the $R$-modules $\Tor_i^R(N,M)$ and $\Ext^i_R(N,M)$, using this
free resolution. \qed\\

\section{Cofiniteness of local cohomology modules}

In this section we give a sufficient condition for a given ideal $I$ of a Noetherian ring $R$ being in $\mathscr{I}(R)$. The main goal of this section is Theorem 3.8, which is a generalization of some results given in \cite{B, BN, Ch, De, DM, Ha, HK, Yo}.

\begin{lem}
\label{3.1}
Let $I$ be an ideal of $R$ and $M$ be a finitely generated $R$-module. Then
$\mathfrak{B}(I,M)\subseteq \big(\mathfrak{C}(I,M)\cup\mathfrak{D}(I,M)\big).$ In particular,
$$\mAss_R M=\mathfrak{A}(I,M) \cup \mathfrak{B}(I,M)\cup \mathfrak{C}(I,M)\cup \mathfrak{D}(I,M)$$ if and only if $\mAss_R M=\mathfrak{A}(I,M) \cup \mathfrak{C}(I,M) \cup \mathfrak{D}(I,M).$
\end{lem}
\proof Let $\p \in \mathfrak{B}(I,M)$. Then $\p\in \mAss_R M$ and $\cd(I,R/\p)=1$. So, we have $\p\not \in V(I)$ and hence $H^0_I(R/\p)=0$. Now, we consider the following two cases:\\

Case 1. Assume that $H^1_I(R/\p)$ is Artinian. Then the $R$-modules $H^i_I(R/\p)$ are Artinian for all $i\geq0$ and hence $q(I,R/\p)=-\infty$. Thus, in view of \cite[Lemma 4.1]{B} we have $\dim R/(I+\p)=0$. Therefore,  $\p\in \mathfrak{D}(I,M)$.\\

Case 2. Assume that $H^1_I(R/\p)$ is not Artinian. Then it is clear that $q(I,R/\p)=1$ and so $\p\in \mathfrak{C}(I,M)$. \qed\\

\begin{lem}
\label{3.2}
Let $I$ be an ideal of $R$. Then the following statements are equivalent:
\begin{enumerate}[\upshape (i)]
  \item  $\mAss_R M=\mathfrak{A}(I,M) \cup \mathfrak{C}(I,M) \cup \mathfrak{D}(I,M)$ for every finitely generated $R$-module $M$.
  \item $\mAss_R R=\mathfrak{A}(I,R)\cup \mathfrak{C}(I,R) \cup \mathfrak{D}(I,R).$
  \item  $\mAss_R M=\mathfrak{A}(I,M) \cup \mathfrak{B}(I,M)\cup \mathfrak{C}(I,M)\cup \mathfrak{D}(I,M)$ for every finitely generated $R$-module $M$.
  \item $\mAss_R R=\mathfrak{A}(I,R) \cup \mathfrak{B}(I,R) \cup \mathfrak{C}(I,R)\cup \mathfrak{D}(I,R).$
 \end{enumerate}
 \end{lem}
\proof (i)$\Rightarrow$(ii) It is trivial.\\

(ii)$\Rightarrow$(i) Let $M$ be a finitely generated $R$-module. If $M=0$ then the assertion is clear. So, we assume that $M\neq 0$. It is enough to prove that
$$\mAss_R M\subseteq\big(\mathfrak{A}(I,M)\cup\mathfrak{C}(I,M)\cup\mathfrak{D}(I,M)\big).$$
Let $\q \in \mAss_R M.$ Then there exists an element $\p\in \mAss_R R$ such that $\q\in V(\p)$. From the hypothesis it follows that $$\p \in \big(\mathfrak{A}(I,R)\cup\mathfrak{C}(I,R)\cup\mathfrak{D}(I,R)\big).$$  We consider the following three cases:\\

Case 1. Assume that $\p \in \mathfrak{A}(I,R)$. If $I+\p=R$ then $R=I+\p\subseteq I+\q$ and hence $I+\q=R$. Therefore, $\q\in \mathfrak{A}(I,M)$. Also, if $\p\supseteq I$ then $I\subseteq \p\subseteq \q$ and hence $\q\in \mathfrak{A}(I,M)$.\\

Case 2. Assume that $\p \in \mathfrak{C}(I,R)$. Then, as $$\Supp R/\q  \subseteq \Supp R/\p,$$ \cite[Theorem 3.2]{DY} implies that   $q(I,R/\q)\leq q(I,R/\p)= 1$. If $q(I,R/\q)=1$, then  $\q\in \mathfrak{C}(I,M)$. Also, if $q(I,R/\q)=0$ then it is clear that $\q\in V(I)$ and hence $\q\in \mathfrak{A}(I,M)$. Also, if $q(I,R/\q)=-\infty$ then in view of \cite[Lemma 4.1]{B} we have $\q\in \mathfrak{D}(I,M)\cup\mathfrak{A}(I,M)$.\\

Case 3. Assume that $\p \in \mathfrak{D}(I,R)$. Then, by the definition we have $$\dim R/(I+\q)\leq \dim R/(I+\p)\leq 1,$$which implies that $\q \in \mathfrak{D}(I,M)\cup\mathfrak{A}(I,M)$.\\

(i)$\Leftrightarrow$(iii) and (ii)$\Leftrightarrow$(iv) Follow from Lemma 3.1. \qed\\

The following well known lemma is needed in the proof of Lemma 3.4.

\begin{lem}
\label{3.3}
Let $(R,\m)$ be a Noetherian complete local ring and $I$ be an ideal of $R$. Assume that $\cd(I,R)=t\geq0$ and suppose that the $R$-module $H^t_I(R)$ is Artinian and $I$-cofinite. Then $$\Att_R H^t_I(R)=\{\q \in \mAss_R R\,\,:\,\,\dim R/\q =t\,\,{\rm and}\,\,\Rad(I+\q)=\m\}.$$
\end{lem}
\proof See \cite[Lemma 2.3]{B0}.\qed\\

\begin{lem}
\label{3.4}
Let $(R,\m)$ be a Noetherian complete local ring, $I$ be an ideal of $R$ and $M$ be a finitely generated $R$-module. Then  $\mathfrak{C}(I,M)\subseteq \big(\mathfrak{B}(I,M)\cup\mathfrak{D}(I,M)\big).$ In particular,
$$\mAss_R M=\mathfrak{A}(I,M) \cup \mathfrak{B}(I,M)\cup \mathfrak{C}(I,M)\cup \mathfrak{D}(I,M)$$ if and only if $\mAss_R M=\mathfrak{A}(I,M) \cup \mathfrak{B}(I,M) \cup \mathfrak{D}(I,M).$
\end{lem}
\proof Let $\p \in \mathfrak{C}(I,M)$. Then $\p\in \mAss_R M$ and $q(I,R/\p)=1$. So, we have $\cd(I,R/\p)\geq q(I,R/\p)=1$. Now, we consider the following two cases:\\

Case 1. Assume that $\cd(I,R/\p)=1$. Then $\p \in \mathfrak{B}(I,M)$.\\

Case 2. Assume that $\cd(I,R/\p)=t>1$. Then by \cite[Theorem 4.9]{B} the $R$-module $H^t_I(R/\p)\simeq H^t_{(I+\p)/\p}(R/\p)$ is Artinian and $(I+\p)/\p$-cofinite. So, Lemma 3.3 yields that $\dim R/\p=t$ and $\Rad(I+\p)=\m$ and hence $\dim R/(I+\p)=0$. Thus $\p\in \mathfrak{D}(I,M)$.\qed\\

The following corollary is a consequence of Lemmas 3.1, 3.2 and 3.4.

\begin{cor}
\label{3.5}
Let $(R,\m)$ be a Noetherian complete local ring and $I$ be an ideal of $R$. Then the following statements are equivalent:
\begin{enumerate}[\upshape (i)]
   \item  $\mAss_R R=\mathfrak{A}(I,R) \cup \mathfrak{C}(I,R) \cup \mathfrak{D}(I,R)$.
  \item $\mAss_R R=\mathfrak{A}(I,R)\cup \mathfrak{B}(I,R) \cup \mathfrak{D}(I,R).$
  \item  $\mAss_R M=\mathfrak{A}(I,M) \cup \mathfrak{B}(I,M) \cup \mathfrak{D}(I,M)$, for each finitely generated $R$-module $M$.
  \item $\mAss_R R=\mathfrak{A}(I,R) \cup \mathfrak{B}(I,R) \cup \mathfrak{C}(I,R)\cup \mathfrak{D}(I,R).$
  \item  $\mAss_R M=\mathfrak{A}(I,M) \cup \mathfrak{B}(I,M) \cup \mathfrak{C}(I,M)\cup \mathfrak{D}(I,M)$, for each finitely generated $R$-module $M$.
 \end{enumerate}
 \end{cor}
\proof (i)$\Leftrightarrow$(iv) Follows from Lemma 3.1.\\

(ii)$\Leftrightarrow$(iv) and (iii)$\Leftrightarrow$(v) Follow from Lemma 3.4.\\

(iv)$\Leftrightarrow$(v) Follows from Lemma 3.2.\qed\\

Combining Corollary 3.5 and Theorem 2.5 we have the following result.

\begin{cor}
 \label{3.6}
Let $(R,\m)$ be a Noetherian complete local ring and $I$ be an ideal of $R$ such that $$\mAss_R R=\mathfrak{A}(I,R)\cup\mathfrak{B}(I,R)\cup\mathfrak{C}(I,R)\cup\mathfrak{D}(I,R).$$  Then
$\mathscr{C}(R, I)_{cof}$ is Abelian.
\end{cor}
\qed\\

The following proposition plays a key role in the proof of the main result of this section.

\begin{prop}
\label{3.7}
Let $R$ be a Noetherian ring and $I$ be an ideal of $R$. Let $M$ be a finitely generated $R$-module such that $$\mAss_R M=\mathfrak{A}(I,M) \cup\mathfrak{C}(I,M)\cup\mathfrak{D}(I,M).$$
    Then the $R$-modules $H^i_I(M)$ are $I$-cofinite for all $i\geq 0$.
\end{prop}
\proof
For each finitely generated $R$-module $N$, set $$\mathfrak{A}_1(I,N):=\{\p \in \mAss_R N\,\,:\,\,I+\p=R\}\,\,\,{\rm and}\,\,\,\mathfrak{A}_2(I,N):=\{\p \in \mAss_R N\,\,:\,\,\p\supseteq I\}.$$ In order to prove the assertion, without loss of generality, we may assume that $\mathfrak{A}_1(I,M)=\emptyset$. Because, in the case where $\mathfrak{A}_1(I,M)\neq \emptyset$, we can see that $I\Gamma_{J_1}(M)=\Gamma_{J_1}(M) $, where $J_1:=\bigcap_{\p \in \mathfrak{A}_1(I,M)}\p$. Hence,  $H^i_I(\Gamma_{J_1}(M))=0$ for each integer $i\geq0$. Consequently, the exact sequence $$0\longrightarrow \Gamma_{J_1}(M) \longrightarrow M \longrightarrow
M/\Gamma_{J_1}(M) \longrightarrow 0,$$induces the isomorphisms $H^i_I(M)\simeq H^i_I(M/\Gamma_{J_1}(M))$ for all integers $i\geq0$. Since $$\Ass_R M/\Gamma_{J_1}(M)=\Ass_R M \backslash V(J_1)$$ and $(\mathfrak{C}(I,M)\cup\mathfrak{D}(I,M))\cap\mathfrak{A}_1(I,M)=\emptyset$, it follows that
\begin{eqnarray*}
\mAss_R M/\Gamma_{J_1}(M)&=&\big(\mathfrak{A}_2(I,M) \cup\mathfrak{C}(I,M)\cup\mathfrak{D}(I,M)\big)\\&=& \big(\mathfrak{A}_2(I,M/\Gamma_{J_1}(M)) \cup\mathfrak{C}(I,M/\Gamma_{J_1}(M))\cup\mathfrak{D}(I,M/\Gamma_{J_1}(M))\big).
\end{eqnarray*}

So, replacing $M$ by $M/\Gamma_{J_1}(M)$, without loss of generality, we may assume $\mathfrak{A}_1(I,M)=\emptyset$.\\

Moreover, again without loss of generality, we may assume that $\mathfrak{A}_2(I,M)=\emptyset$. Because, in the case where $\mathfrak{A}_2(I,M)\neq \emptyset$, we can see that $\Gamma_{J_2}(M)\subseteq\Gamma_{I}(M)$, where $J_2:=\cap_{\p \in \mathfrak{A}_2(I,M)}\p$. Therefore, $H^i_I(\Gamma_{J_2}(M))=0$ for each integer $i\geq1$.  Consequently, the exact sequence $$0\longrightarrow \Gamma_{J_2}(M) \longrightarrow M \longrightarrow
M/\Gamma_{J_2}(M) \longrightarrow 0,$$yields the isomorphisms $H^i_I(M)\simeq H^i_I(M/\Gamma_{J_2}(M))$ for all integers $i\geq1$. Since, $$\Ass_R M/\Gamma_{J_2}(M)=\Ass_R M \backslash V(J_2)$$ it follows that $$\mAss_R M/\Gamma_{J_2}(M)=\big(\mathfrak{C}(I,M)\cup\mathfrak{D}(I,M)\big)\backslash V(J_2)=\big(\mathfrak{C}(I,M/\Gamma_{J_2}(M))\cup\mathfrak{D}(I,M/\Gamma_{J_2}(M))\big).$$
Furthermore, it is clear that the finitely generated $I$-torsion module $H^0_I(M)$ is $I$-cofinite. So, replacing $M$ by $M/\Gamma_{J_2}(M)$, without loss of generality, we may assume that $\mathfrak{A}_2(I,M)=\emptyset.$ Therefore, using the fact that $$\mathfrak{A}(I,M)=\mathfrak{A}_1(I,M)\cup \mathfrak{A}_2(I,M),$$ without loss of generality we may assume $\mathfrak{A}(I,M)=\emptyset.$ \\

Next, let $M$ be a finitely generated $R$-module with $$\mAss_R M=\big(\mathfrak{C}(I,M)\cup\mathfrak{D}(I,M)\big)\neq \emptyset.$$

Henceforth, we shall prove the assertion for all possible cases. To do this, we consider the following three cases:\\

Case 1. Assume that $\mathfrak{C}(I,M)\neq\emptyset$ and $\mathfrak{D}(I,M)=\emptyset$. Then $T:=\bigoplus_{\p\in \mAss_R M}R/\p$ is a finitely generated
 $R$-module with $\Supp T = \Supp M$. So, it follows from \cite[Theorem 3.2]{DY}
  that $$q(I,M)=q(I,T)={\rm max}\{q(I,R/\p)\,\,:\,\,\p\in \mathfrak{C}(I,M)\}=1.$$ Thus, the assertion follows from
   \cite[Theorem 4.9]{B}.\\

Case 2. Assume that $\mathfrak{C}(I,M)=\emptyset$ and $\mathfrak{D}(I,M)\neq\emptyset$. Then for each $\q \in \Supp M/IM$ there exists $\p \in \mAss_R M=\mathfrak{D}(I,M)$ such that $\q \in \big(V(\p)\cap V(I)\big)=V(I+\p)$ and so $$\dim R/\q\leq \dim R/(I+\p)\leq 1.$$ Hence $\dim M/IM = {\rm sup}\{\dim R/\q\,\,:\,\,\q \in \Supp M/IM\}\leq 1.$ So, the assertion follows from \cite[Corollary 2.7]{BN}. \\

Case 3. Assume that $\mathfrak{C}(I,M)\neq\emptyset$ and $\mathfrak{D}(I,M)\neq\emptyset$. Set $$J_3:=\bigcap_{\p\in \mathfrak{C}(I,M)}\p,\,\,\,\,\,K:=\Gamma_{J_3}(M)\,\,\,\,\,{\rm and}\,\,\,\,\,L:=\bigoplus_{\p\in \mathfrak{C}(I,M)}R/\p.$$ Then, it follows from \cite[Theorem 3.2]{DY}
  that $$q(I,L)={\rm max}\{q(I,R/\p)\,\,:\,\,\p\in \mathfrak{C}(I,M)\}= 1.$$ Moreover, since $\Supp K \subseteq V(J_3)=\Supp L$, it follows from
  \cite[Theorem 3.2]{DY} that $$q(I,K)\leq  q(I,L)= 1.$$ Hence, in view of  \cite[Theorem 4.9]{B} the $R$-module $H^i_I(K)$ is Artinian and $I$-cofinite, for each $i\geq 2$. Furthermore, the exact sequence $$0\longrightarrow K \longrightarrow M \longrightarrow
M/K \longrightarrow 0,$$induces an exact sequence $$H^i_{I}(K)\longrightarrow H^i_{I}(M)\longrightarrow H^i_{I}(M/K)\longrightarrow H^{i+1}_{I}(K),$$
for each integer $i\geq 2$. Therefore, applying \cite[Corollary 4.4]{Me}, it follows that the $R$-modules $H^i_{I}(M)$ are $I$-cofinite for all integers $i\geq 2$, if and only if, the $R$-modules $H^i_{I}(M/K)$ are $I$-cofinite for all integers $i\geq 2$.\\

On the other hand, $$\Ass_R M/K=\Ass_R M/\Gamma_{J_3}(M)=\Ass_R M \backslash V(J_3)$$ and so that $\mAss_R M/K=\mathfrak{D}(I,M)\backslash V(J_3)=\mathfrak{D}(I,M/K)$. Thus, by the proof of Case 2, the $R$-modules $H^i_I(M/K)$ are $I$-cofinite for all integers $i\geq 2$. Now, we are ready to deduce that the $R$-modules $H^i_{I}(M)$ are $I$-cofinite for all integers $i\geq 2$.  Also, since the $R$-module $H^0_I(M)$
  is finitely generated with support in $V(I)$ it follows that $H^0_I(M)$ is $I$-cofinite. Therefore, for each $i\neq 1$ the $R$-module
  $H^i_I(M)$ is $I$-cofinite. Hence, by \cite[Proposition 3.11]{Me} the $R$-module $H^1_I(M)$ is $I$-cofinite too.\qed\\

Now, we are ready to deduce the main result of this section.

\begin{thm}
\label{3.8}
Let $I$ be an ideal of $R$ such that $$\mAss_R R=\mathfrak{A}(I,R)\cup \mathfrak{B}(I,R)\cup\mathfrak{C}(I,R)\cup\mathfrak{D}(I,R).$$ Then $I\in\mathscr{I}(R)$.
 \end{thm}
\proof The assertion follows from Lemma 3.2 and Proposition 3.7.\qed\\

\section{Cofiniteness and cohomological dimension}

In this section we give a formula for the cohomological dimensions of finitely generated modules over a Noetherian complete local ring $R$ with respect to ideals in $\mathscr{I}(R)$. The main goal of this section is Theorem 4.5. The following lemmas will be quite useful in this section.

\begin{lem}
\label{4.1}
 Let $(R,\m)$ be a Noetherian local ring  and $I$ be an ideal of $R$. Let $M$ be a non-zero $I$-cofinite $R$-module of dimension $n$. Then $\cd(\m,M)=n$.
 \end{lem}
\proof Follows from \cite[Theorem 2.9]{M} and the {\it Grothendieck's Vanishing Theorem}.\qed\\

\begin{lem}
\label{4.2}
 Let $(R,\m)$ be a Noetherian local ring  and $I$ be an ideal of $R$. Assume that $M$ is a non-zero $I$-cofinite $R$-module of dimension $n\geq 1$ and $x\in \m$ is an element with the property that $\Ass_R M \cap V(Rx)\subseteq \{\m\}$. Then the $R$-module $H^0_{Rx}(M)\simeq H^0_{I+Rx}(M)$ is Artinian and $I$-cofinite and $\dim H^1_{Rx}(M)=n-1$. Moreover, the $R$-modules $H^0_{Rx}(M)$ and $H^1_{Rx}(M)\simeq H^1_{I+Rx}(M)$ are $(I+Rx)$-cofinite. In particular, $H^1_{Rx}(M)\neq0$.
\end{lem}
\proof
 Since by the hypothesis $M$ is $I$-cofinite, from the definition it follows that $M$ is an $I$-torsion $R$-module. Thus, by \cite[Exercise 2.1.9]{BS} we have that $H^i_{Rx}(M)\simeq H^i_{I+Rx}(M)$ for each $i\geq 0$. From the hypothesis $\Ass_R M \cap V(Rx)\subseteq \{\m\}$, it follows that $\Supp H^0_{Rx}(M)\subseteq \{\m\}$. Now as $(0:_{H^0_{Rx}(M)}I)\subseteq (0:_MI)$ and by the hypothesis the $R$-module $$(0:_MI)\simeq \Hom_R(R/I,M)$$ is finitely generated, it follows that the $R$-module $(0:_{H^0_{Rx}(M)}I)$ is of finite length. Therefore, the result \cite[Proposition 4.1]{Me} implies that the $R$-module $\Gamma_{Rx}(M)\simeq H^0_{Rx}(M)$ is Artinian and $I$-cofinite. In particular, $\dim M/\Gamma_{Rx}(M)=\dim M=n \geq 1$. Also, another usage of \cite[Proposition 4.1]{Me} yields that the $R$-module $\Gamma_{Rx}(M)\simeq H^0_{Rx}(M)$ is $(I+Rx)$-cofinite. Now the exact sequence $$0\longrightarrow \Gamma_{Rx}(M) \longrightarrow M \longrightarrow
M/\Gamma_{Rx}(M) \longrightarrow 0,$$ yields that the $R$-module $M/\Gamma_{Rx}(M)$ is $I$-cofinite too. Therefore, considering the relation $I \subseteq (I+Rx)$, it follows from \cite[Corollary 1]{DM} or \cite[Corollary 2.5]{Me}, that for all $i \geq 0$,
the $R$-modules ${\rm Ext}^{i}_{R}(R/(I+Rx),M/\Gamma_{Rx}(M))$ are finitely generated. By \cite[Remark 2.2.17]{BS}, there is an exact sequence $$0 \longrightarrow M/\Gamma_{Rx}(M) \longrightarrow M_x \longrightarrow H^1_{Rx}(M)\longrightarrow 0.\,\,\,\,\,(4.2.1)$$ On the other hand,
 multiplication by $x$ is an automorphism on $M_x$.  Therefore, multiplication by $x$ is an automorphism on
  $\Ext_R^i(R/(I+Rx),M_x)$, for all $i\geq 0$. But, since $x\in (I+Rx)$ it follows that, multiplication by $x$ on $\Ext_R^i(R/(I+Rx),M_x)$
  is the zero map, for all $i\geq 0$. Thus, $\Ext_R^i(R/(I+Rx),M_x)=0$ for all $i\geq 0$. Furthermore, for each integer $i\geq0$, the exact sequence $(4.2.1)$ induces an exact sequence $$\Ext_R^i(R/(I+Rx),M_x) \longrightarrow \Ext_R^i(R/(I+Rx),H^1_{Rx}(M)) \longrightarrow$$$$ \Ext_R^{i+1}(R/(I+Rx),M/\Gamma_{Rx}(M))  \longrightarrow \Ext_R^{i+1}(R/(I+Rx),M_x),$$which yields the isomorphisms $$\Ext_R^i(R/(I+Rx),H^1_{Rx}(M)) \simeq \Ext_R^{i+1}(R/(I+Rx),M/\Gamma_{Rx}(M)),$$ for all $i\geq 0$. This means that the $(I+Rx)$-torsion $R$-module $H^1_{Rx}(M)\simeq H^1_{I+Rx}(M)$ is $(I+Rx)$-cofinite. \\
 Also,  since for all $i\geq 0$ the multiplication by $x$ is an automorphism on
  $H^i_{\m}(M_x)$ and the $R$-module $H^i_{\m}(M_x)$ is $Rx$-torsion, it follows that $H^i_{\m}(M_x)=0$ for all $i\geq 0$. Moreover, for each integer $i\geq0$, the exact sequence $(4.2.1)$ induces an exact sequence $$H^i_{\m}(M_x) \longrightarrow H^i_{\m}(H^1_{Rx}(M)) \longrightarrow H^{i+1}_{\m}(M/\Gamma_{Rx}(M))  \longrightarrow H^{i+1}_{\m}(M_x),$$which yields the isomorphisms $$H^i_{\m}(H^1_{Rx}(M))\simeq H^{i+1}_{\m}(M/\Gamma_{Rx}(M)),$$ for all $i\geq 0$. Therefore, applying Lemma 4.1 it follows that
\begin{eqnarray*}
\dim H^1_{Rx}(M)&=&\cd(\m,H^1_{Rx}(M))\\&=&\cd(\m,M/\Gamma_{Rx}(M))-1\\&=&\dim M/\Gamma_{Rx}(M)-1\\&=&\dim M-1\\&=&n-1.
\end{eqnarray*}
Now, the proof is complete.\qed\\

The following proposition and its corollary play a key role in the proof of Theorem 4.5.

\begin{prop}
\label{4.3}
  Let $(R,\m)$ be a Noetherian complete local domain  and $I$ be a proper ideal of $R$. Assume that the $R$-module $H^i_I(R)$ is $I$-cofinite for each $i\geq0$. Then $\cd(I,R)=\height I$.
\end{prop}
\proof In view of {\it Grothendieck's Vanishing Theorem} we have $\cd(I,R)\leq \dim R$ and so $\dim R-\cd(I,R)\geq 0$. Now in order to prove the assertion we use induction on $n:=\dim R-\cd(I,R).$ If $n=0$, then $H^{\dim R}_{I}(R)\neq 0$, and hence by the {\it Lichtenbaum-Hartshorne Vanishing Theorem} we have $\Rad(I)=\m$. Thus, $\height I=\height \m=\dim R=\cd(I,R)$. Suppose,
 inductively, that $0<n \leq \dim R$ and the result has been proved for
 $n-1$. Then from the {\it Lichtenbaum-Hartshorne Vanishing Theorem} it follows that $\dim R/I> 0$, and hence $$\m\not\subseteq \bigcup_{\frak P\in \Assh_R R/I}\frak P.$$ Also, since for all $i\geq 0$ the $R$-modules $H^i_I(R)$ are $I$-cofinite, it follows that the set $$T:=\bigcup_{i\geq 0}\Ass_R H^i_I(R)= \bigcup_{i= 0}^{\dim R}\Ass_R H^i_I(R)$$is finite. Therefore, there exists an element $x\in \m$ such that $$x\not\in \left(\bigcup_{\Q\in (T\backslash\{\m\})}\Q\right)\bigcup\left(\bigcup_{\frak P\in \Assh_R R/I}\frak P\right).$$ Then, in view of the Lemma 4.2, the $R$-modules $$H^0_{Rx}(H^i_I(R))\simeq H^0_{I+Rx}(H^i_I(R))\,\,\,\,{\rm and}\,\,\,\ H^1_{Rx}(H^i_I(R))\simeq H^1_{I+Rx}(H^i_I(R))$$ are $(I+Rx)$-cofinite for all $i\geq 0$. Furthermore, by \cite[Corollary 3.5]{Sc} for each $i\geq 0$, there exists an exact sequence $$0\longrightarrow H^1_{Rx} (H^{i-1}_I(R)) \longrightarrow H^{i}_{I+Rx}(R)\longrightarrow H^0_{Rx} (H^{i}_I(R))\longrightarrow 0.\,\,\,\,\,(4.3.1)$$These exact sequences together with Lemma 4.2 imply that the $R$-modules $H^{i}_{I+Rx}(R)$ are $(I+Rx)$-cofinite for all $i\geq 0$. Now, we claim that $\cd(I+Rx,R)=\cd(I,R)+1$. Considering the exact sequences $(4.3.1)$ for all $i\geq 0$ and Lemma 4.2 it is enough to prove $\dim H^{\cd(I,R)}_I(R)\geq 1$. Assume the opposite. Then by \cite[Proposition 4.1]{Me} the non-zero $R$-module $H^{\cd(I,R)}_I(R)$ is Artinian and $I$-cofinite. Since, by the hypothesis $R$ is a domain, it follows from Lemma 3.3 that $\Rad(I)=\m$, which is a contradiction. Therefore, $\cd(I+Rx,R)=\cd(I,R)+1$ and so by the inductive hypothesis we have $$\height(I+Rx)=\cd(I+Rx,R)=\cd(I,R)+1.$$ Since $R$ is a catenary domain, it follows that $$\dim R/(I+Rx)=\dim R -\height(I+Rx)=\dim R -\cd(I+Rx,R).$$ But, we have $$x\not\in \bigcup_{\frak P\in \Assh_R R/I}\frak P.$$Therefore, $$\dim R- \height I=\dim R/I= \dim R/(I+Rx)+1=\dim R-\cd(I,R),$$ which implies that $$\height I=\cd(I,R).$$ This completes the inductive step. \qed\\

Note that if $I\in\mathscr{I}(R)$ then it follows from the {\it Independence Theorem} and Lemma 2.4 that $(I+J)/J\in\mathscr{I}(R/J)$, for every ideal $J$ of $R$. Henceforth, we shall use this fact several times.

\begin{cor}
\label{4.4}
 Let $(R,\m)$ be a Noetherian complete local ring  and $I\in\mathscr{I}(R)$.  Then $\cd(I,R/\p)=\height (I+\p)/\p$, for each $\p \in \Spec R$.
\end{cor}
\proof Let $\p \in \Spec R$. From the hypothesis $I\in\mathscr{I}(R)$ it follows that $(I+\p)/\p\in\mathscr{I}(R/\p)$. So, the assertion follows from the Proposition 4.3 using the fact that $\cd(I,R/\p)=\cd((I+\p)/\p,R/\p)$.\qed\\

Now, we are ready to state and prove the main result of this section.

\begin{thm}
\label{4.5}
 Let $(R,\m)$ be a Noetherian complete local ring  and $I\in\mathscr{I}(R)$.  Then,
  $$\cd(I,M)={\rm max}\{\height (I+\p)/\p\,\,:\,\,\p\in \mAss_R M\},$$for every non-zero finitely generated $R$-module $M$.
 \end{thm}
\proof
Set $T:=\bigoplus_{\p\in \mAss_R M}R/\p$. Then $T$ is a finitely generated
 $R$-module with $\Supp T = \Supp M$. So, it follows from \cite[Theorem 2.2]{DNT} and Corollary 4.4
  that
\begin{eqnarray*}
\cd(I,M)&=&\cd(I,T)\\&=&{\rm max}\{\cd(I,R/\p)\,\,:\,\,\p\in \mAss_R M\}\\&=&{\rm max}\{\height (I+\p)/\p\,\,:\,\,\p\in \mAss_R M\}.
\end{eqnarray*}
\qed\\

Recall that if $R$ is a Noetherian ring of finite Krull dimension, then we say that $R$ is  {\it equidimensional}
if $\mAss_R R=\Assh_R R$.

\begin{prop}
\label{4.6}
 Let $(R,\m)$ be a Noetherian complete local ring  and $I\in\mathscr{I}(R)$.  Then for each $\p\in \Spec R$, the quotient ring $R/(\p+I)$ is equidimensional.
 \end{prop}
\proof
Let $\p \in \Spec R$. Then, $\cd(I,R/\p)=\height (I+\p)/\p$, by Corollary 4.4. Now, if $\mAss_R R/(I+\p)\neq\Assh_R R/(I+\p)$, then there is an element $\q\in \mAss_R R/(I+\p)$ such that $\q\not \in \Assh_R R/(I+\p)$. So $\dim R/\q< \dim R/(I+\p)$, and hence using the fact that $R/\p$ is a catenary domain it follows that $\height \q/\p>\height (I+\p)/\p$. But,  {\it Grothendieck's Non-vanishing Theorem} yields the inequality $$\height (I+\p)/\p=\cd(I,R/\p)\geq\height \q/\p,$$ which is a contradiction.\qed\\

Using an example given in \cite{Ha}, we can construct an example of Noetherian complete local domain $(R,\m)$ of dimension 4, such that $R$ has an ideal $I$ with $\height I=2$ and $I\not\in\mathscr{I}(R)$. Maybe the same property holds in general for any ideal of height $2$ in any Noetherian complete local domain of dimension 4;  because there is no evidence to reject it. Now consider the following two questions: \\

{\bf Question A:} {\it Let $(R,\m)$ be a Noetherian complete local ring and $I\in\mathscr{I}(R)$. Whether $\mAss_R R= \mathfrak{A}(I,R) \cup\mathfrak{B}(I,R) \cup\mathfrak{D}(I,R)$?} \\

{\bf Question B:} {\it Let $(R,\m)$ be a Noetherian complete local domain of dimension $4$ and $I$ be an ideal of $R$ with $\height I=2$. Whether $I\not\in\mathscr{I}(R)$?} \\

 \begin{prop}
 \label{4.7}
Question {\rm A} has an affirmative answer in general if and only if Question {\rm B} has so.
\end{prop}

\proof
$''\Rightarrow''$. Let $(R,\m)$ be a Noetherian complete local domain of dimension $4$ and $I$ be an ideal of $R$ with $\height I=2$. Since, $\height I=2$ it follows from the {\it Grothendieck's Non-vanishing Theorem} that $\cd(I,R)\geq 2$. Therefore, $\{0\}\not \in \big(\mathfrak{A}(I,R) \cup\mathfrak{B}(I,R)\big)$. Also, using the fact that $R$ is a catenary domain it follows that $\dim R/I= \dim R - \height I=4-2=2,$ which implies that $\{0\}\not \in \mathfrak{D}(I,R)$.  So that $$\mAss_R R=\{0\}\not\subseteq \big(\mathfrak{A}(I,R) \cup\mathfrak{B}(I,R) \cup\mathfrak{D}(I,R)\big).$$ Now, it follows from the hypothesis that $I\not\in\mathscr{I}(R)$.\\

 $''\Leftarrow''$. Let $(R,\m)$ be a Noetherian complete local ring and $I\in\mathscr{I}(R)$. Then, by using induction on $d=\dim R$ we prove that $$\mAss_R R= \mathfrak{A}(I,R) \cup\mathfrak{B}(I,R) \cup\mathfrak{D}(I,R).$$ For $d=0,1,2$, the assertion follows from {\it Grothendieck's Vanishing Theorem} and {\it Lichtenbaum-Hartshorne Vanishing Theorem}.  Now assume that $d=3$. In order to prove the assertion it is enough to prove $$\left(\mAss_R R \backslash\big(\mathfrak{A}(I,R) \cup\mathfrak{B}(I,R)\big)\right)\subseteq\mathfrak{D}(I,R).$$ Let $$\p \in \left(\mAss_R R\backslash\big(\mathfrak{A}(I,R) \cup\mathfrak{B}(I,R)\big)\right).$$ Then it follows from  {\it Grothendieck's Vanishing Theorem} that $2\leq \cd(I,R/\p)\leq 3$. If $\cd(I,R/\p)=3$ then the {\it Lichtenbaum-Hartshorne Vanishing Theorem} yields that $\dim R/(I+\p)=0$ and so $\p\in \mathfrak{D}(I,R)$. Now assume that $\cd(I,R/\p)=2$. Then, Corollary 4.4 yields that $\height (I+\p)/\p = \cd(I,R/\p)=2$ and hence $$\dim R/(I+\p)=\dim R/\p - \height (I+\p)/\p\leq 1.$$ Thus, $\p\in \mathfrak{D}(I,R)$.\\

Now assume that $d=4$. In order to prove the assertion it is enough to prove $$\left(\mAss_R R \backslash\big(\mathfrak{A}(I,R) \cup\mathfrak{B}(I,R)\big)\right)\subseteq\mathfrak{D}(I,R).$$ Let $$\p \in \left(\mAss_R R\backslash\big(\mathfrak{A}(I,R) \cup\mathfrak{B}(I,R)\big)\right).$$ Since, $(I+\p)/\p\in\mathscr{I}(R/\p)$, considering the previous lines of the proof, without loss of generality we may assume that $\dim R/\p=4$.  Then it follows from the {\it Grothendieck's Vanishing Theorem} that $2\leq \cd(I,R/\p)\leq 4$. If $\cd(I,R/\p)=4$ then the {\it Lichtenbaum-Hartshorne Vanishing Theorem} yields that $\dim R/(I+\p)=0$ and so $\p\in \mathfrak{D}(I,R)$. Now assume that $\cd(I,R/\p)<4$. Then, we claim that $\cd(I,R/\p)=3$. Assume the opposite. Then, $\cd(I,R/\p)=2$ and so by Corollary 4.4, $\height (I+\p)/\p = \cd(I,R/\p)=2$. But,  $(I+\p)/\p\in\mathscr{I}(R/\p)$ and $\dim R/\p=4$, which is a contradiction.  So, we have $\cd(I,R/\p)=3$ and hence Corollary 4.4 yields that $\height (I+\p)/\p = \cd(I,R/\p)=3$. Therefore, $$\dim R/(I+\p)=\dim R/\p - \height (I+\p)/\p= 1.$$ Thus, $\p\in \mathfrak{D}(I,R)$.\\

Now suppose, inductively, that $d>4$ and the result has been proved for all smaller values of $d$. Then it is enough to prove  $$\left(\mAss_R R \backslash\big(\mathfrak{A}(I,R) \cup\mathfrak{B}(I,R)\big)\right)\subseteq\mathfrak{D}(I,R).$$ Let $$\p \in \left(\mAss_R R\backslash\big(\mathfrak{A}(I,R) \cup\mathfrak{B}(I,R)\big)\right).$$ We claim that $\p\in \mathfrak{D}(I,R).$ Assume the opposite. Then, as $(I+\p)/\p\in\mathscr{I}(R/\p)$ it follows from the inductive hypothesis that $\dim R/\p=d$. Set $t:= \cd(I,R/\p)$. Then, it follows from Corollary 4.4 that $\height (I+\p)/\p=\cd(I,R/\p)=t>1$. In particular, from the fact that $R/\p$ is a catenary domain of dimension $d$ we have that $d-t=\dim R/(I+\p)>1$.\\

Pick  $\q_1\in \Spec R$ with $$\p\subset \q_1 \subset \bigcup_{\Q\in \Assh_R R/(I+\p)}\Q\,\,\,\,{\rm and}\,\,\,\,\height \q_1/\p=1.$$  Since, $\height (I+\q_1)/\p=\height (I+\p)/\p=t$ and $R/\p$ is a catenary domain it follows from the hypothesis that $\height(I+\q_1)/\q_1=t-1$. Then, in view of Collorary 4.4 we have $$\cd(I, R/\q_1)=\height (I+\q_1)/\q_1=t-1>0.$$ Since, $\dim R/(I+\q_1)=\dim R/(I+\p)$ and by the hypothesis
$\p\not\in \mathfrak{D}(I,R)$ it follows that $$\dim R/(I+\q_1)=\dim R/(I+\p)>1.$$ Since by the hypothesis $I\in\mathscr{I}(R)$ it follows that  $(I+\q_1)/\q_1\in\mathscr{I}(R/\q_1)$. Now, applying the inductive hypothesis for the Noetherian complete local domain $R/\q_1$ of dimension $d-1$, it follows that, $$t-1=\cd(I, R/\q_1)=\cd((I+\q_1)/\q_1,R/\q_1)=1,$$ and so $\height (I+\p)/\p=t=2$.\\

Next, let $\frak P\in V(I+\p)$ be a prime ideal with $\height \frak P/(I+\p)=1.$ As $R/\p$ is a catenary domain, using Proposition 4.6, it follows that
\begin{eqnarray*}
\dim R/\frak P&=&\dim R/\p - \height \frak P/\p\\&=&\dim R/\p - (\height (I+\p)/\p+\height \frak P/(I+\p))\\&=&d-3.
\end{eqnarray*}
Pick an element $x\in \frak P$ with $x\not\in \bigcup_{\Q\in \Assh_R R/(I+\p)}\Q.$ Then $\frak P$ contains a prime ideal $\q_2\in \mAss_R R/(\p+Rx)$ with $\height \q_2/\p=1$.
Since, $x\not\in \bigcup_{\Q\in \Assh_R R/(I+\p)}\Q$ and $x\in \m$ it follows that $$\dim R/(I+\p+Rx)= \dim R/(I+\p)-1=d-3=\dim R/\frak P.$$ Also, it is clear that $$I+\p+Rx\subseteq I+\q_2\subseteq \frak P$$ and hence
$$d-3=\dim R/\frak P\leq \dim R/(I+\q_2)\leq \dim R/(I+\p+Rx)=d-3,$$which implies that $$\dim R/\frak P= \dim R/(I+\q_2)=d-3.$$ Thus, using the fact that $R/\p$ is a catenary domain we get
$$\height (I+\q_2)/\p=\dim R/\p - \dim R/(I+\q_2)=3.$$
As $R/\p$ is a catenary domain, using Proposition 4.6, it follows that $$\height (I+\q_2)/\q_2=\height (I+\q_2)/\p- \height \q_2/\p=3-1=2.$$ Then, Corollary 4.4 yields that $$\cd(I, R/\q_2)=\height (I+\q_2)/\q_2=2.$$ On the other hand, since by the hypothesis $I\in\mathscr{I}(R)$ it follows that $(I+\q_2)/\q_2\in\mathscr{I}(R/\q_2)$. Considering the hypothesis $$\cd((I+\q_2)/\q_2, R/\q_2)=\cd(I, R/\q_2)=2>1$$  and  applying the inductive hypothesis for the Noetherian complete local domain $R/\q_2$ of dimension $d-1$ it follows that $$d-3=(d-1)-2=\dim R/\q_2-\height (I+\q_2)/\q_2=\dim R/(I+\q_2)\leq 1.$$ Whence, we have $d\leq4$. Now, we have achieved
the desired contradiction. This completes the inductive step. \qed\\

We close this section by the following three questions.\\

{\bf Question C:} {\it Whether $\mAss_R R= \mathfrak{A}(I,R) \cup\mathfrak{B}(I,R)\cup\mathfrak{C}(I,R)  \cup\mathfrak{D}(I,R)$, for each $I\in\mathscr{I}(R)$?} \\

{\bf Question D:} {\it Let $I$ be an ideal of $R$ with $$\mAss_R R= \mathfrak{A}(I,R) \cup\mathfrak{B}(I,R)\cup\mathfrak{C}(I,R)  \cup\mathfrak{D}(I,R).$$ Whether the category $\mathscr{C}(R, I)_{cof}$ is Abelian?}\\

{\bf Question E:} {\it Let $I$ be an ideal of $R$. Whether the category  $\mathscr{C}(R, I)_{cof}$ is Abelian if and only if $I\in\mathscr{I}(R)$?}\\

\subsection*{Acknowledgements}
The author would like to acknowledge his deep gratitude from
the referee for a very careful reading of the manuscript
and many valuable suggestions. He also, would like
to thank to School of Mathematics, Institute for Research in Fundamental
Sciences (IPM) for its financial support.


\end{document}